\documentclass[12pt]{article}
\usepackage{amssymb}
\begin{document}

%Macro for the adjoint:
\def\ad{\mbox{ad}\,}
%Macro for the black-board bold
\def\b#1{{\mathbb #1}}
%Macro for the calligraphic letters
\def\c#1{{\cal #1}}
%Macro to enumerate formulae according to sections and subsections
\newcommand{\sect}[1]{\setcounter{equation}{0}\section{#1}}
\renewcommand{\theequation}{\thesection.\arabic{equation}}
\newcommand{\subsect}[1]{\setcounter{equation}{0}\subsection{#1}}
\renewcommand{\theequation}{\thesection.\arabic{equation}}
\def\1{{\bf 1}}
\newcommand{\be}{\begin{equation}}
\newcommand{\ee}{\end{equation}}
\newcommand{\ba}{\begin{array}}
\newcommand{\ea}{\end{array}}

\title{Leibniz Rules and Reality Conditions}

\author{Gaetano Fiore$^1$, \  John Madore$^{2,3}$
        \and
        $\strut^1$Dip. di Matematica e Applicazioni, Fac.  di Ingegneria\\ 
        Universit\`a di Napoli, V. Claudio 21, 80125 Napoli
        \and
        $\strut^2$Max-Planck-Institut f\"ur Physik 
        (Werner-Heisenberg-Institut)\\
        F\"ohringer Ring 6, D-80805 M\"unchen
        \and
        $\strut^3$Laboratoire de Physique Th\'eorique et Hautes Energies\\
        Universit\'e de Paris-Sud, B\^atiment 211, F-91405 Orsay
        }

\date{}

\maketitle

\abstract{An analysis is made of reality conditions within the context
of noncommutative geometry. We show that if a covariant derivative
satisfies a given left Leibniz rule then a right Leibniz rule is
equivalent to the reality condition. We show also that the matrix
which determines the reality condition must satisfy the Yang-Baxter
condition if the extension of the covariant derivative to tensor
products is to satisfy the reality condition. This is equivalent to the
braid condition for the matrix which determines the right Leibniz rule.}

\vfill
\noindent
Preprint 98-13, Dip. Matematica e Applicazioni, Universit\`a di Napoli\\
\medskip
\eject

\parskip 4pt plus2pt minus2pt

\sect{Introduction and motivation}

In noncommutative geometry (or algebra), reality conditions are not as
natural as they can be in the commutative case; the product of two
hermitian elements is no longer necessarily hermitian. The product of
two hermitian differential forms is also not necessarily hermitian. It is
our purpose here to analyze this problem in some detail. If the reality
condition is to be extended to a covariant derivative then we shall show
that there is a unique correspondence between its existence and the
existence of a left and right Leibniz rule.  We shall show also that the
matrix which determines the reality condition must satisfy the
Yang-Baxter condition if the extension of the covariant derivative to
tensor products is to be well-defined. This is equivalent to the braid
condition for the matrix which determines the right Leibniz rule. It is
necessary in discussing the reality of the curvature form.

There is not as yet a completely satisfactory definition of either a
linear connection or a metric within the context of noncommutative
geometry but there are definitions which seem to work in certain cases.
In the present article we chose one particular
definition~\cite{DubMadMasMou96}. We refer to a recent review
article~\cite{Mad97} for a list of some other examples and references to
alternative definitions. More details of one alternative version can be
found, for example, in the book by Landi~\cite{Lan97}. For a general
introduction to more mathematical aspects of the subject we refer to the
book by Connes~\cite{Con94}.  Although we expect our results to have a
more general validity we shall prove them only in a particular version
of noncommutative geometry which can be considered as a noncommutative
extension of the moving-frame formalism of E. Cartan.  This implies that
we suppose that the module of 1-forms is free as a right or left
module. As a bimodule it will always be projective with one generator,
the generalized `Dirac operator'.  More details can be found
elsewhere~\cite{Mad95, DubMadMasMou96}.  We shall use here the
expression `connection' and `covariant derivative' synonymously.

In the second section we describe briefly what we mean by the frame
formalism and we recall the particular definition of a covariant
derivative which we use. In the third section we discuss the reality
condition. We describe here the relation between the map which
determines the right Leibniz rule and the map which determines the
reality condition. The last section contains a generalization to
higher tensor powers.

\sect{The frame formalism}

The starting point is a noncommutative algebra $\c{A}$ and over
$\c{A}$ a differential calculus~\cite{Con94} $\Omega^*(\c{A})$. 
We recall that a differential calculus is completely
determined by the left and right module structure of the
$\c{A}$-module of 1-forms $\Omega^1(\c{A})$.  We shall restrict our
attention to the case where this module is free of rank $n$ as a left
or right module and possesses a special basis $\theta^a$, 
$1\leq a \leq n$, which commutes with the elements $f$ of the algebra:
\be
[f, \theta^a] = 0.                                             \label{fund} 
\ee
In particular, if the geometry has a commutative limit then the
associated manifold must be parallelizable.  We shall refer to the
$\theta^a$ as a `frame' or `Stehbein'.  The integer $n$ plays the role
of `dimension'; it can be greater than the dimension of the limit
manifold but in this case the frame will have a singular limit. We
suppose further~\cite{DimMad96} that the basis is dual to a set of inner
derivations $e_a = \ad \lambda_a$. This means that the differential is
given by the expression
\be
df = e_a f \theta^a = [\lambda_a, f] \theta^a.                 \label{defdiff}
\ee
One can rewrite this equation as 
\be
df = -[\theta,f],                                              \label{extra}
\ee
if one introduces~\cite{Con94} the `Dirac operator'
\be
\theta = - \lambda_a \theta^a.                                 \label{dirac}
\ee

There is a bimodule map $\pi$ of the space 
$\Omega^1(\c{A}) \otimes_\c{A} \Omega^1(\c{A})$ onto the space
$\Omega^2(\c{A})$ of 2-forms and we can write
\be
\theta^a \theta^b = P^{ab}{}_{cd} \theta^c \otimes \theta^d     \label{proj}
\ee
where, because of (\ref{fund}), the $P^{ab}{}_{cd}$ belong to the center
$\c{Z}(\c{A})$ of $\c{A}$. We shall suppose that the center is
trivial, $\c{Z}(\c{A}) = \b{C}$, and therefore the
components $P^{ab}{}_{cd}$ are complex numbers. Define the Maurer-Cartan 
elements $C^a{}_{bc} \in \c{A}$ by the equation
\be
d\theta^a = - {1\over 2} C^a{}_{bc} \theta^b \theta^c.
\ee
Because of~(\ref{proj}) we can suppose that 
$C^a{}_{bc} P^{bc}{}_{de} = C^a{}_{de}$. It follows from the equation
$d(\theta^a f - f \theta^a) = 0$ that there exist elements
$F^a{}_{bc}$ of the center such that
\be
C^a{}_{bc} = F^a{}_{bc} - 2 \lambda_e P^{(ae)}{}_{bc}        \label{consis1}
\ee
where $(ab)$ means symmetrization of the indices $a$ and $b$.
If on the other hand we define $K_{ab}$ by the equation
\be
d\theta + \theta^2 = {1\over 2} K_{ab} \theta^a \theta^b,      \label{consis2}
\ee
then if follows from (2.3) and the identity $d^2 = 0$ that the $K_{ab}$
must belong to the center. Finally it can be shown~\cite{DimMad96,
MadMou98} that in order that~(\ref{consis1}) and (\ref{consis2}) be
consistent with one another the original $\lambda_a$ must satisfy the
condition
\be
2 \lambda_c \lambda_d P^{cd}{}_{ab} - 
\lambda_c F^c{}_{ab} - K_{ab} = 0.                             \label{manca}
\ee
This gives to the set of $\lambda_a$ the structure of a twisted Lie
algebra with a central extension.

We propose as definition of a linear connection a map~\cite{Kos60,
CunQui95}
\be
\Omega^1(\c{A}) \buildrel D \over \longrightarrow 
\Omega^1(\c{A}) \otimes_\c{A} \Omega^1(\c{A})                  \label{2.2.4} 
\ee 
which satisfies both a left Leibniz rule
\be
D (f \xi) =  df \otimes \xi + f D\xi                           \label{2.2.2}
\ee
and a right Leibniz rule~\cite{DubMadMasMou96}
\be
D(\xi f) = \sigma (\xi \otimes df) + (D\xi) f                 \label{second}
\ee
for arbitrary $f \in \c{A}$ and $\xi \in \Omega^1(\c{A})$.  
We have here introduced a generalized permutation
\be
\Omega^1(\c{A}) \otimes_\c{A} \Omega^1(\c{A})
\buildrel \sigma \over \longrightarrow
\Omega^1(\c{A}) \otimes_\c{A} \Omega^1(\c{A})                  \label{2.2.5}
\ee
in order to define a right Leibniz rule which is consistent with the
left one, It is necessarily bilinear.  A linear connection is therefore
a couple $(D, \sigma)$. It can be shown that a necessary as well as
sufficient condition for torsion to be right-linear is that $\sigma$
satisfy the consistency condition
\be
\pi \circ (\sigma + 1) = 0.                                    \label{2.2.6}
\ee
Using the fact that $\pi$ is a projection one sees that the most general
solution to this equation is given by
\be
1 + \sigma = ( 1 - \pi) \circ \tau                             \label{2.2.7}
\ee
where $\tau$ is an arbitrary bilinear map
\be
\Omega^1(\c{A}) \otimes \Omega^1(\c{A})
\buildrel \tau \over \longrightarrow
\Omega^1(\c{A}) \otimes \Omega^1(\c{A}).                         \label{2.2.8}
\ee
If we choose $\tau = 2$ then we find $\sigma = 1 - 2 \pi$ and 
$\sigma^2 = 1$. The eigenvalues of $\sigma$ are then equal to $\pm 1$.  
The map~(\ref{2.2.4}) has a natural extension~\cite{Kos60}
\be
\Omega^*(\c{A}) \buildrel D \over \longrightarrow 
\Omega^*(\c{A}) \otimes_\c{A} \Omega^1(\c{A})                  \label{2.2.4ex}
\ee 
to the entire tensor algebra given by a graded Leibniz rule.

This general formalism can be applied in particular to differential
calculi with a frame.  Since $\Omega^1(\c{A})$ is a free module the
maps $\sigma$ and $\tau$ can be defined by their action on the basis
elements:
\be
\sigma (\theta^a \otimes \theta^b) = 
S^{ab}{}_{cd} \theta^c \otimes \theta^d,   \qquad
\tau (\theta^a \otimes \theta^b) = 
T^{ab}{}_{cd} \theta^c \otimes \theta^d.                         \label{2.2.9}
\ee
By the sequence of identities
\be
f S^{ab}{}_{cd} \theta^c \otimes \theta^d = 
\sigma (f \theta^a \otimes \theta^b) = 
\sigma (\theta^a \otimes \theta^b f) =
S^{ab}{}_{cd} f \theta^c \otimes \theta^d                        \label{2.2.10}
\ee
and the corresponding ones for $T^{ab}{}_{cd}$ we conclude that the
coefficients $S^{ab}{}_{cd}$ and $T^{ab}{}_{cd}$ must lie in 
$\c{Z}(\c{A})$. From~(\ref{2.2.7}) the most general form for 
$S^{ab}{}_{cd}$ is
\be
S^{ab}{}_{cd} = 
T^{ab}{}_{ef} (\delta^e_c \delta^f_d - P^{ef}{}_{cd}) 
- \delta^a_c \delta^b_d.                                         \label{2.2.11}
\ee

A covariant derivative can be defined also by its action on the basis 
elements:
\be
D\theta^a = - \omega^a{}_{bc} \theta^b \otimes \theta^c.       \label{2.2.12}
\ee
The coefficients here are elements of the algebra. They are restricted
by (2.1) and the the two Leibniz rules. The torsion 2-form is defined as
usual as
\be
\Theta^a = d \theta^a - \pi \circ D \theta^a.             \label{deftorsion}
\ee
If $F^a{}_{bc} = 0$ then it is easy to check~\cite{DubMadMasMou96} 
that
\be
D_{(0)} \theta^a = - \theta \otimes \theta^a + 
\sigma (\theta^a \otimes \theta)                            \label{2.2.14}
\ee
defines a torsion-free covariant derivative.                    
The most general $D$ for fixed $\sigma$ is of the form
\be
D = D_{(0)} + \chi                                           \label{2.2.15}
\ee
where $\chi$ is an arbitrary bimodule morphism
\be
\Omega^1(\c{A}) \buildrel \chi \over \longrightarrow
\Omega^1(\c{A}) \otimes \Omega^1(\c{A}).               \label{2.2.16}
\ee
If we write
\be
\chi (\theta^a) = - \chi^a{}_{bc} \theta^b \otimes \theta^c   \label{2.2.17}
\ee
we conclude that $\chi^a{}_{bc} \in \c{Z}(\c{A})$. 
In general a covariant derivative is torsion-free provided the condition
\be
\omega^a{}_{de} P^{de}{}_{bc} = {1\over 2}C^a{}_{bc}          \label{2.2.19}
\ee
is satisfied. The covariant derivative~(\ref{2.2.15}) is torsion free 
if and only
if
\be
\pi \circ \chi = 0.                                            \label{2.2.20}
\ee

One can define a metric by the condition
\be
g(\theta^a \otimes \theta^b) = g^{ab}                           \label{2.2.21}
\ee
where the coefficients $g^{ab}$ are elements of $\c{A}$.  To be well
defined on all elements of the tensor product $\Omega^1(\c{A})
\otimes_\c{A} \Omega^1(\c{A})$ the metric must be bilinear and by
the sequence of identities
\be
f g^{ab} = g(f \theta^a \otimes \theta^b) 
= g(\theta^a \otimes \theta^b f) = g^{ab} f                     \label{2.2.22}
\ee
one concludes that the coefficients must lie in $\c{Z}(\c{A})$.
We define the metric to be symmetric if
\be 
g \circ \sigma \propto g.
\ee
This is a natural generalization of the situation in ordinary
differential geometry where symmetry is respect to the flip which
defines the forms.  If $g^{ab} = g^{ba}$ then by a linear transformation
of the original $\lambda_a$ one can make $g^{ab}$ the components of the
Euclidean (or Minkowski) metric in dimension $n$. It will not
necessarily then be symmetric in the sense that we have just used the
word.

The covariant derivative~(\ref{2.2.12}) is compatible with the metric if
and only if~\cite{DimMad96}
\be
\omega^a{}_{bc} + \omega_{cd}{}^e S^{ad}{}_{be} = 0.            \label{2.2.23}
\ee
This is a `twisted' form of the usual condition that
$g_{ad}\omega^d{}_{bc}$ be antisymmetric in the two indices $a$ and $c$
which in turn expresses the fact that for fixed $b$ the
$\omega^a{}_{bc}$ form a representation of the Lie algebra of the
Euclidean group $SO(n)$ (or the Lorentz group). 
 When $F^a{}_{bc} = 0$
the condition that~(\ref{2.2.12}) be metric compatible can be
written~\cite{DimMad96} as
\be
S^{ae}{}_{df} g^{fg} S^{bc}{}_{eg} = g^{ab} \delta^c_d.         \label{2.2.24}
\ee

Introduce the standard notation $\sigma_{12} = \sigma \otimes 1$,
$\sigma_{23} = 1 \otimes \sigma$, to extend to three factors of a module
any operator $\sigma$ defined on a tensor product of two factors.
Then there is a natural continuation of the map (\ref{2.2.4}) to the 
tensor product $\Omega^1(\c{A}) \otimes_\c{A} \Omega^1(\c{A})$ 
given by the map
\be
D_2(\xi \otimes \eta) = D\xi \otimes \eta + 
\sigma_{12} (\xi \otimes D\eta)                                 \label{2.2.4e}
\ee
The map $D_2 \circ D$ has no nice properties but if one introduces the
notation $\pi_{12} = \pi \otimes 1$ then by analogy with the commutative 
case one can set
\be
D^2 = \pi_{12} \circ D_2 \circ D
\ee
and formally define the curvature as the map
\be
\mbox{Curv}:\, \Omega^1(\c{A}) \longrightarrow 
\Omega^2(\c{A}) \otimes_\c{A} \Omega^1(\c{A})                   \label{curv}
\ee
given by $\mbox{Curv} = D^2$.  This coincides with the composition of
the first two maps of the series of~(\ref{2.2.4ex}).  Because of the
condition (\ref{2.2.6}) Curv is left linear. It can be written out in
terms of the frame as
\be
\mbox{Curv} (\theta^a) =
- {1 \over 2} R^a{}_{bcd} \theta^c \theta^d \otimes \theta^b   \label{2.16}
\ee
Similarly one can define a Ricci map
\be
\mbox{Ric} (\theta^a) =
{1 \over 2} R^a{}_{bcd} \theta^c g(\theta^d \otimes \theta^b).  \label{2.17}
\ee
It is given by
\be
\mbox{Ric} \, (\theta^a) = R^a{}_b \theta^b.                   \label{2.18}
\ee
The above definition of curvature is not satisfactory in the
noncommutative case~\cite{DubMadMasMou96}. For example, from~(\ref{2.16})
one sees that Curv can only be right linear if 
$R^a{}_{bcd} \in \c{Z}(\c{A})$

The curvature $\mbox{Curv}_{(0)}$ of the covariant derivative $D_{(0)}$
defined in (\ref{2.2.14}) can be readily calculated. One finds after a
short calculation that it is given by the expression
\be
\mbox{Curv}_{(0)} (\theta^a) = \theta^2 \otimes \theta^a +
\pi_{12} \sigma_{12}\sigma_{23} \sigma_{12}
(\theta^a \otimes \theta \otimes \theta).                 \label{2.19}
\ee
If  $\xi = \xi_a \theta^a$ is a general 1-form then since Curv is left
linear one can write
\be
\mbox{Curv}_{(0)} (\xi ) =  \xi_a \theta^2 \otimes \theta^a +
\pi_{12} \sigma_{12}\sigma_{23} \sigma_{12}
(\xi \otimes \theta \otimes \theta).                       \label{2.20}
\ee
The lack of right-linearity of Curv is particularly evident in
this last formula.

\sect{The involution}

Suppose now that $\c{A}$ is a $*$-algebra. We would like to choose the
differential calculus such that the reality condition $(df)^* = df^*$
holds. This can at times be difficult~\cite{olezu}. We must
require that the derivations $e_a$ satisfy the reality condition
\be
(e_a f^*)^* = e_a f
\ee
which in turn implies that the $\lambda_a$ are antihermitian. One finds
that for general $f \in \c{A}$ and $\xi \in \Omega^1(\c{A})$ one has
\be
(f \xi)^* = \xi^* f^*, \qquad (\xi f)^* = f^* \xi^*.
\ee
 From the duality condition we find that
\be
(\theta^a)^* = \theta^a, \qquad \theta^* = - \theta.      \label{RealTheta}
\ee
There are elements $I^{ab}{}_{cd}, J^{ab}{}_{cd} \in \c{Z}(\c{A})$
such that
\be
(\theta^a \theta^b)^* = \imath(\theta^a \theta^b) = 
I^{ab}{}_{cd} \theta^c \theta^d, \qquad
(\theta^a \otimes \theta^b)^* = \jmath_2(\theta^a \otimes \theta^b) = 
J^{ab}{}_{cd} \theta^c \otimes \theta^d.
\ee
We can suppose that
\be
I^{ab}{}_{cd} P^{cd}{}_{ef} = I^{ab}{}_{ef}.
\ee
We have then
\be
(I^{ab}{}_{cd})^* I^{cd}{}_{ef} = P^{ab}{}_{ef}, \qquad
(J^{ab}{}_{cd})^* J^{cd}{}_{ef} = \delta^a_e \delta^b_f.
\ee
The compatibility condition with the product
\be
\pi \circ \jmath_2 = \imath \circ \pi 
\ee
becomes
\be
(P^{ab}{}_{cd})^* J^{cd}{}_{ef} =  I^{ab}{}_{cd} P^{cd}{}_{ef} 
 = I^{ab}{}_{ef}.             \label{compatibility}
\ee
But since the frame is hermitian and associated to derivations we have 
from~(\ref{defdiff})
\be
(dfdg)^*  = [d(f\,dg)]^* = d(f\,dg)^* = d(dg^*\,f) = - dg^*df^*
\ee
for arbitrary $f$ and $g$. It follows that
\be 
(e_a f e_b g)^* I^{ab}{}_{cd} \theta^c \theta^d = 
- e_b g^* e_a f^* \theta^b \theta^a = 
- (e_a f e_b g)^* \theta^b \theta^a
\ee
and we must conclude that 
\be
I^{ab}{}_{cd} = - P^{ba}{}_{cd}.                                \label{I-P}
\ee
It can be shown~\cite{MadMou98} that the right-hand satisfies a weak
form of the Yang-Baxter equation, which would imply some sort of braid
condition on the left-hand side.

The compatibility condition with the product implies then that
\be
(P^{ab}{}_{cd})^* P^{dc}{}_{ef} = P^{ba}{}_{ef}.
\ee
For general $\xi, \eta \in \Omega^1(\c{A})$ it follows from~(\ref{I-P})
that
\be
(\xi \eta)^* = - \eta^* \xi^*.                               \label{antisym}
\ee
In particular
\be
(\theta^a \theta^b)^* = - \theta^b \theta^a. 
\ee
The product of two frame elements is hermitian then if and only if they
anticommute.  More generally one can extent the involution to the entire
algebra of forms by setting
\be
(\alpha \beta)^* =
(-1)^{pq} \beta^* \alpha^*                                \label{sign}
\ee
if $\alpha \in \Omega^p(\c{A}$) and $\beta \in \Omega^q(\c{A})$.
When the frame exists one has necessarily also the relations
\be
(f \xi \eta)^* = (\xi \eta)^* f^*, \qquad
(f \xi \otimes \eta)^* = (\xi \otimes \eta)^* f^*
\ee
for arbitrary $f \in \c{A}$.  If in particular $P^{ab}{}_{cd}$ is given
by
\be
P^{ac}{}_{cd} = 
{1 \over 2} (\delta^a_c \delta^b_d - \delta^a_d \delta^b_c) 
\ee
we can choose $\jmath_2$ to be the identity. In this case the
$F^c{}_{ab}$ are hermitian and the $K_{ab}$ anti-hermitian elements of 
$\c{Z}(\c{A})$. An involution can be introduced on the algebra of
forms even if they are not defined using derivations~\cite{Con95}

We require that the metric be real; if $\xi$ and $\eta$ are hermitian
1-forms then $g(\xi \otimes \eta)$ should be an hermitian element of
the algebra. The  reality condition for the metric becomes therefore
\be
g((\xi \otimes \eta)^*) = 
(g(\xi \otimes \eta))^*                                      \label{reality}
\ee 
and puts further constraints
\be
S^{ab}{}_{cd} g^{cd} = (g^{ba})^*
\ee
on the matrix of coefficients $g^{ab}$.

We shall also require the reality condition
\be
D\xi^* = (D\xi)^*                                          \label{basic}
\ee
on the connection,
which can be rewritten also in the form
\be
D\circ\jmath_1=\jmath_2\circ D.                            \label{prima}
\ee
%%
%%so as to start to give the example of the lowest $n$...
This must be consistent with the Leibniz rules.
There is little one can conclude in general but if the differential is
based on real derivations then from the equalities 
\be
(D(f\xi))^* = D((f\xi)^*) = D(\xi^* f^*)
\ee
one finds the conditions
\be
(df \otimes \xi)^* + (f D\xi)^* = 
\sigma(\xi^* \otimes df^*) + (D\xi^*) f^*.
\ee 
Since this must be true for arbitrary $f$ and $\xi$ we conclude that
\be
(df \otimes \xi)^* = \sigma(\xi^* \otimes df^*)
\ee
and
\be
(f D\xi)^* = (D\xi^*) f^*.
\ee
We shall suppose~\cite{DubMadMasMou95, KasMadTes97} that the
involution is such that in general
\be
(\xi \otimes \eta)^* = \sigma(\eta^* \otimes \xi^*).          \label{TPI}
\ee
A change in $\sigma$ therefore implies a change in the definition of
an hermitian tensor. From the compatibility
conditions~(\ref{compatibility}) and~(\ref{2.2.6}) one can
deduce~(\ref{antisym}).  The condition that the star operation be in
fact an involution places a constraint on the map $\sigma$:
\be
(\sigma(\eta^* \otimes \xi^*))^* =  (\xi \otimes \eta).    \label{invconst}
\ee

It is clear that there is an intimate connection between the reality
condition and the right-Leibniz rule.
The expression~(\ref{TPI}) for the involution on tensor products becomes
the identity
\be
J^{ab}{}_{cd} = S^{ba}{}_{cd}.                                  \label{J-S} 
\ee
This is consistent with~(\ref{I-P}) because of~(\ref{2.2.6}). It forces also
the constraint
\be
(S^{ba}{}_{cd})^* S^{dc}{}_{ef} = \delta^a_e \delta^b_f      \label{s-unitary}
\ee
on $\sigma$. Equation~(\ref{J-S}) can be also read from right to left as a 
definition of the right-Leibniz rule in terms of the hermitian structure.

The condition that the connection~(\ref{2.2.12}) be real can be
written as
\be
(\omega^a{}_{bc})^* = \omega^a{}_{de} (J^{de}{}_{bc})^*.      \label{real1st}
\ee
One verifies immediately that the connection~(\ref{2.2.14})
is real. 

In order for the curvature to be real we must require that the extension 
of the involution to the tensor product of three elements of
$\Omega^1(\c{A})$ be such that 
\be
\pi_{12} \circ D_2(\xi \otimes \eta)^* =
\Big(\pi_{12} \circ D_2(\xi \otimes \eta)\Big)^*.
\ee
We shall impose a stronger condition. We shall require that $D_2$ be real:
\be
D_2(\xi \otimes \eta)^* = (D_2(\xi \otimes \eta))^*.          \label{strong}
\ee
This condition can be made more explicit when a frame exists. In this
case the map $D_2$ is given by
\be
D_2 (\theta^a \otimes \theta^b) = 
- (\omega^a{}_{pq} \delta^b_r + S^{ac}{}_{pq} \omega^b{}_{cr})
\theta^p \otimes \theta^q \otimes \theta^r.
\label{explici}
\ee
To solve the reality condition~(\ref{strong})
we introduce elements $J^{abc}{}_{def} \in \c{Z}(\c{A})$ such that
\be
(\theta^a \otimes \theta^b \otimes \theta^c)^* = 
\jmath_3(\theta^a \otimes \theta^b \otimes \theta^c) =
J^{abc}{}_{def} \theta^d \otimes \theta^e \otimes \theta^f. \label{involution}
\ee
Using~(\ref{J-S}) one finds then that the equality
\be
D_2 \circ \jmath_2 = \jmath_3 \circ D_2
\label{dopo}
\ee
can be written in the form
\be
J^{ab}{}_{pq} (\omega^p{}_{de} \delta^q_f + J^{rp}{}_{de} \omega^q{}_{rf}) =
\Big((\omega^a{}_{pq})^* \delta^b_r + 
(J^{sa}{}_{pq})^* (\omega^b{}_{sr})^*\Big) J^{pqr}{}_{def}.   \label{4.2.43}
\ee
This equation must be solved for $J^{abc}{}_{def}$ as a function of
$J^{ab}{}_{cd}$. One cannot simply cancel the factor  
$\omega^a{}_{bc}$ since it satisfies constraints. As a test case we 
choose~(\ref{2.2.14}). We find that~(\ref{4.2.43}) is satisfied
provided 
\be
J^{abc}{}_{def} = J^{ab}{}_{pq}J^{pc}{}_{dr}J^{qr}{}_{ef}
                = J^{bc}{}_{pq}J^{aq}{}_{rf}J^{rp}{}_{de}.     \label{Y-B}
\ee
The second equality is the Yang-Baxter Equation written out 
with indices. Using this equation it follows that~(\ref{involution}) is 
indeed an involution:
\be
(J^{abc}{}_{pqr})^* J^{pqr}{}_{def} = \delta^a_d \delta^b_e \delta^c_f.
\ee
Using Equation~(\ref{Y-B}) the Equation~(\ref{4.2.43}) can be written
in the form
\be
J^{ab}{}_{pe} \omega^p{}_{cd} - J^{ap}{}_{de} \omega^b{}_{cp} +
J^{ab}{}_{pq} J^{rp}{}_{cd} \omega^q{}_{re} - 
J^{qb}{}_{cp} J^{rp}{}_{de} \omega^a{}_{qr} = 0.               \label{real2nd}
\ee
The connection then must satisfy two reality conditions,
Equation~(\ref{real1st}) and Equation~(\ref{real2nd}).  The second
condition can be rewritten more concisely in the form
\be
D_2 \circ \sigma = \sigma_{23} \circ D_2.                     \label{equi}
\ee
In fact, using Equations~(\ref{explici}), (\ref{2.2.2}) one finds
\[
\begin{array}{l}
D_2\Big(\sigma(f\theta^b\otimes\theta^a)\Big) - 
\sigma_{23} \circ \Big(D_2(f\theta^b\otimes\theta^a)\Big) =\\
S^{ba}{}_{pq}D_2(f\theta^p\otimes\theta^q)- df \otimes 
\sigma(\theta^b\otimes\theta^a)-
f(\omega^b{}_{cp} \delta^a_r + S^{bq}{}_{cp} \omega^a{}_{qr})
\sigma_{23}(\theta^c \otimes \theta^p \otimes \theta^r) = \\
f\Big(S^{ba}{}_{pq}(\omega^p{}_{cd} \delta^q_e + S^{pr}{}_{cd} 
\omega^q{}_{re}) - (\omega^b{}_{cp} \delta^a_r + S^{bq}{}_{cp} 
\omega^a{}_{qr}) S^{pr}{}_{de}\Big)\theta^c\otimes\theta^d\otimes\theta^e.
\end{array}
\]
Because of (\ref{J-S}), the right-hand side of this equation vanishes
if and only if the left-hand side of Equation~(\ref{real2nd}) is zero.
One can check that equations (\ref{equi}) and (\ref{strong}) are
equivalent, once the definitions of $\jmath_2,\jmath_3$ and the 
property~(\ref{basic}) are postulated.

It is reasonable to suppose that even in the absence of a frame the
constraints~(\ref{s-unitary}) and the Yang-Baxter condition hold. 
The former has in fact already been written~(\ref{invconst}) in general.
The map $\jmath_3$ can be written as
\be
(\xi\otimes\eta\otimes\zeta)^* \equiv
\jmath_3(\xi\otimes\eta\otimes\zeta)
= \sigma_{12}\sigma_{23} \sigma_{12}
(\zeta^*\otimes\eta^*\otimes\xi^*).                           \label{defj3}
\ee
Because of~(\ref{J-S}) the Yang-Baxter condition for $\jmath_2$ becomes
the braid equation
\be
\sigma_{12}\sigma_{23}\sigma_{12}=\sigma_{23}\sigma_{12}
\sigma_{23}                                                  \label{genbraid}
\ee
for the map $\sigma$.

\sect{Higher tensor and wedge powers}

Just as we have~(\ref{2.2.4e}) defined $D_2$ we can 
introduce a set $D_n$ of covariant derivatives
\be
D_n: \, \bigotimes_1^n \Omega^1(\c{A}) \longrightarrow
\bigotimes_1^{n+1} \Omega^1(\c{A}) 
\ee
for arbitrary integer $n$ by using $\sigma$ to place the operator $D$ in
its natural position to the left. For instance,
\be
D_3 = \Big(D\otimes 1 \otimes 1 + 
\sigma_{12}(1 \otimes D\otimes 1)+
\sigma_{12}\sigma_{23}(1\otimes 1\otimes D)\Big)             \label{example}
\ee
If the condition~(\ref{genbraid}) is satisfied then these $D_n$ will
also be real in the sense that
\be
D_n \circ \jmath_n = \jmath_{n+1} \circ D_n              \label{genreality1}
\ee
where the $\jmath_n$ are the natural extensions of $\jmath_2$ and
$\jmath_3$. For instance, $\jmath_4$ is defined by
\be
(\xi\otimes\eta\otimes\zeta\otimes \omega)^* \equiv
\jmath_4(\xi\otimes\eta\otimes\zeta\otimes \omega) =
\sigma_{12}\sigma_{23}\sigma_{12}\sigma_{34}\sigma_{23}\sigma_{12}
(\omega^*\otimes\zeta^*\otimes\eta^*\otimes\xi^*).               \label{defj4}
\ee
The general rule to construct $\jmath_n$ is the following. Let
$\epsilon$ denote the ``flip'', the permutator of two objects,
$\epsilon(\xi\otimes\eta)=\eta\otimes\xi$, and more generally let
$\epsilon_n$ denote the inverse-order permutator of $n$ objects.  For
instance, the action of $\epsilon_3$ is given by
\be
\epsilon_3(\zeta\otimes\eta\otimes\xi)=\xi\otimes\eta\otimes\zeta.
                                                                \label{deco}
\ee 
The maps $\epsilon,\epsilon_n$ are $\b{C}$-bilinear but not
$\c{A}$-bilinear, and are involutive.  One can decompose $\epsilon_n$
as a product of $\epsilon_{i(i\!+\!1)}$. One finds for $n=3$
\be
\epsilon_3=\epsilon_{12}\epsilon_{23}\epsilon_{12}=
\epsilon_{23}\epsilon_{12}\epsilon_{23}.
\ee
The second equality expresses the fact that $\epsilon$ fulfils the
braid equation.  In a more abstract but compact notation the
definitions (\ref{TPI}), (\ref{defj3}) and (\ref{defj4}) can be written
in the form
\begin{eqnarray}
\jmath_2 &=&\sigma\,\ell_2,                                 \label{defj2abs}\\
\jmath_3 &=&\sigma_{12}\sigma_{23}\sigma_{12} \,\ell_3,     \label{defj3abs}\\
\jmath_4 &=&
\sigma_{12}\sigma_{23}\sigma_{12}\sigma_{34}\sigma_{23}\sigma_{12}
\,\ell_4.
                                                             \label{defj4abs}
\end{eqnarray}
We have here defined the involution on the 1-forms as $\jmath_1$,
and 
\be
\ell_n = (\underbrace{\jmath_1 \otimes \ldots \otimes
\jmath_1}_{\mbox{$n$~times}})\, \epsilon_n.
\ee
The $\ell_n$ is clearly an involution,
since $\epsilon_n$ commutes with the tensor product of the $\jmath_1$'s.  
The products of $\sigma$'s
appearing in the definitions of $\jmath_3,\jmath_4$ are obtained from
the decompositions of $\epsilon_3,\epsilon_4$ by replacing each
$\epsilon_{i(i\!+\!1)}$ by $\sigma_{i(i\!+\!1)}$. In this way,
$\jmath_3,\jmath_4$ have the correct classical limit, since in this
limit $\sigma$ become the ordinary flip $\epsilon$.  In the same way
as different equivalent decompositions of $\epsilon_3,\epsilon_4$ are
possible, so different products of $\sigma$ factors in
(\ref{defj3abs}), (\ref{defj4abs}) are allowed; they are all equal,
once Equation~(\ref{genbraid}) is fulfilled.  The same rules described
for $n=3,4$ should be used to define $\jmath_n$ for $n>4$.
 
The definition of $\jmath_n$ can be given also some equivalent
recursive form which will be useful for the proofs below, namely
\begin{eqnarray}
\jmath_3 &=&
\sigma_{12}\sigma_{23} \epsilon_{23} \epsilon_{12}
(\jmath_1\otimes \jmath_2),                            \label{poi1} \\
\jmath_4 &=&\sigma_{12}\sigma_{23}\sigma_{34} \epsilon_{34} 
\epsilon_{23} \epsilon_{12}(\jmath_1\otimes \jmath_3), \label{poi2} \\
&=&\sigma_{23}\sigma_{34}\sigma_{12}\sigma_{23} \epsilon_{23}
 \epsilon_{12} \epsilon_{34} \epsilon_{23}(\jmath_2\otimes \jmath_2),
                                                        \label{poi3}
\end{eqnarray}
and so forth to higher orders. Again, these definitions are unambiguous
because of the braid equation~(\ref{genbraid}).

Now we wish to show that, if the braid equation is fulfilled and
$\jmath_2$ is an involution, that is, Equation~(\ref{invconst}) is
satisfied then $\jmath_n$ is also an involution for $n>2$. Note that
the constraint~(\ref{invconst}) in the more abstract notation
introduced above becomes
\be
\jmath_2=\jmath_2^{-1}=
\epsilon\circ(\jmath_1\otimes\jmath_1)\circ\sigma^{-1}.    \label{defj2ab}
\ee
As a first step one checks that for $i=1,...,n\!-\!1$
\be
\sigma_{i(i\!+\!1)}\,\ell_n =
\ell_n\,\sigma_{(n\!-\!i)(n\!+\!1\!-\!i)}^{-1}.                \label{fifa}
\ee
The latter relation can be proved recursively. We show in particular
how from the relation with $n=2$ follows the relation with $n=3$:
\be
\begin{array}{lcl}
\sigma_{12}\,\ell_3 
&\stackrel{(\ref{deco})}{=}
&\sigma_{12}(\jmath_1\otimes\jmath_1\otimes\jmath_1)
\,  \epsilon_{12} \epsilon_{23} \epsilon_{12}\,\sigma_{23}\sigma_{23}^{-1}\\
&\stackrel{(\ref{defj2abs})}{=} &
(\jmath_2\otimes\jmath_1)\, \epsilon_{23} \epsilon_{12}\,
\sigma_{23}\sigma_{23}^{-1} \\
&=& (\jmath_2\sigma\otimes\jmath_1)\, 
\epsilon_{23} \epsilon_{12}\,\sigma_{23}^{-1}\\
&\stackrel{(\ref{defj2ab})}{=}
& (\jmath_1\otimes\jmath_1\otimes\jmath_1)
 \epsilon_{12} \epsilon_{23} \epsilon_{12}\,\sigma_{23}^{-1}\\
&\stackrel{(\ref{deco})}{=}& \ell_3\sigma_{23}^{-1}.
\nonumber
\end{array}                                                \label{semplice}
\ee
Now it is immediate to show that $\jmath_n$ is an involution. Again, we
explicitly reconsider the case $n=3$:
\begin{eqnarray}
(\jmath_3)^2 &=&
\sigma_{12}\sigma_{23}\sigma_{12}\,\ell_3\,\sigma_{23}\sigma_{12}
\sigma_{23}\,\ell_3                                               \nonumber\\
&\stackrel{(\ref{semplice})}{=} & \ell_3\,
\sigma^{-1}_{23}\sigma^{-1}_{12}\sigma^{-1}_{23}
\sigma_{23}\sigma_{12}\sigma_{23}\,\ell_3                         \nonumber\\
&=&1.                                                              \nonumber
\end{eqnarray}

In order to prove~(\ref{genreality1}) it is useful to prove first
a direct consequence of relation~(\ref{equi}):
\be
D_n\circ\sigma_{(i\!-\!1)i}=\sigma_{i(i\!+\!1)}\circ D_n.
\label{lemma}
\ee
The recursive proof is straightforward. For instance,
\[
D_3\sigma_{23}=[D\otimes 1 \otimes 1 +\sigma_{12} (1\otimes D_2)]
\sigma_{23}\stackrel{(\ref{equi})}{=}\sigma_{34}(D\otimes 1 \otimes 1) +
\sigma_{12}\sigma_{34}(1\otimes D_2)=\sigma_{34} D_3.
\]
Now (\ref{genreality1}) can be proved recursively. For instance,
\be
\ba{lcl}
D_3 \jmath_3 &\stackrel{(\ref{poi1})}{=}& 
D_3\sigma_{12}\sigma_{23}  \epsilon_{23} \epsilon_{12} 
(\jmath_1\otimes\jmath_2)\\
&\stackrel{(\ref{lemma})}{=}
&\sigma_{23}\sigma_{34} D_3 \epsilon_{23} \epsilon_{12} 
(\jmath_1\otimes\jmath_2)\\
&\stackrel{(\ref{example})}{=}&\sigma_{23}
\sigma_{34} [D_2\otimes 1+\sigma_{12}\sigma_{23}(1\otimes 
1\otimes D)]  \epsilon_{23} \epsilon_{12} 
(\jmath_1\otimes\jmath_2)\\
&=&\sigma_{23}\sigma_{34} [ \epsilon_{34} \epsilon_{23} \epsilon_{12} 
(1\otimes D_2) +\sigma_{12}\sigma_{23} \epsilon_{23} 
\epsilon_{12} \epsilon_{34} \epsilon_{23} (D\otimes 1\otimes 1)] 
(\jmath_1\otimes\jmath_2)\\
&\stackrel{(\ref{dopo})}{=}&
\sigma_{23}\sigma_{34} [ \epsilon_{34} \epsilon_{23} \epsilon_{12} 
(\jmath_1\otimes\jmath_3 D_2) +\sigma_{12}\sigma_{23} \epsilon_{23}
 \epsilon_{12} \epsilon_{34} \epsilon_{23}(\jmath_2 D\otimes\jmath_2)]\\
&\stackrel{(\ref{poi2}),(\ref{poi3})}{=}
&\sigma_{12}^{-1}\jmath_4(1\otimes D_2) +
\jmath_4(D\otimes 1\otimes 1)\\
&=&\jmath_4[\sigma_{12}(1\otimes D_2) +
(D\otimes 1\otimes 1)]\\
&\stackrel{(\ref{example})}{=}& \jmath_4 D_3.
                                                                 \nonumber
\ea
\ee

For the second-last equality we have used the relation
$\sigma_{12}^{-1}\jmath_4 = \jmath_4\sigma_{12}$, which can be easily
proven using Equations~(\ref{genbraid}) and~(\ref{fifa}).

For further developments it is convenient to interpret $\sigma$ as a
``braiding'', in the sense of Majid \cite{majid}. This is possible
because of Eqution~(\ref{genbraid}). In that framework, the bilinear
map $\sigma$ can be naturally extended first to higher tensor powers
of $\Omega^1(\c{A})$,
\be
\sigma:(\underbrace{\Omega^1\otimes\ldots\otimes\Omega^1}_{\mbox{$p$ times}})
\otimes(\underbrace{\Omega^1\otimes\ldots\otimes\Omega^1}_{\mbox{$k$ times}})
\rightarrow\underbrace{\Omega^1\otimes\ldots\otimes\Omega^1}_{\mbox{$p\!+\!k$ 
times}}.
\ee
This extension can be found by applying iteratively the rules
\be
\ba{l}
\sigma\Big((\xi\otimes\eta)\otimes \zeta\Big) =
\sigma_{12}\sigma_{23}(\xi\otimes\eta\otimes \zeta), \\
\sigma\Big(\xi\otimes(\eta\otimes \zeta)\Big)=
\sigma_{23}\sigma_{12}(\xi\otimes\eta\otimes \zeta).
\ea
\ee
Here $\xi,\eta,\zeta$ are elements of three arbitrary tensor powers
of $\Omega^1(\c{A})$.  It is easy to show that there is no ambiguity
in the iterated definitions, and that the extended map still satisfies
the braid equation (\ref{genbraid}).  These are general properties of
a braiding.

Thereafter, by applying $p\!+\!k\!-\!2$ times the projector $\pi$ to
the previous equation, so as to transform the relevant tensor products
into wedge products, $\sigma$ can be extended also as a map
\be
\sigma:\Omega^p(\c{A})\otimes\Omega^k(\c{A})\rightarrow
\Omega^k(\c{A})\otimes \Omega^p(\c{A}).
\ee
For instance, we shall define $\sigma$ on $\Omega^2\otimes\Omega^1$
and $\Omega^1\otimes\Omega^2$ respectively through
\be
\ba{lcl}
\sigma(\xi\eta\otimes \zeta) & =
&\pi_{23}\sigma\Big((\xi\otimes\eta)\otimes \zeta\Big),\\
\sigma(\xi\otimes\eta\zeta) & =
&\pi_{12}\sigma\Big(\xi\otimes(\eta\otimes \zeta)\Big).
\ea
\ee
Under suitable assumptions on $\pi$, the extended $\sigma$ still
satisfies the braid equation (\ref{genbraid}). It follows that the
same formulae presented above in this section can be used to extend
the involutions $\jmath_n$ to tensor powers of higher degree forms in
a compatible way with the action of $\pi$, that is, in such a way that
$\jmath_2\circ\pi_{12}=\pi_{12}\circ\jmath_3$, and so forth. Finally,
also the covariant derivatives $D_n$ can be extended to tensor powers
of higher degree forms in such a way that~(\ref{genreality1}) is still
satisfied. These results will be shown in detail elsewhere.

\section*{Acknowledgment} 
One of the authors (JM) would like to thank J. Wess for his hospitality
and the Max-Planck-Institut f\"ur Physik in M\"unchen for financial
support.

\end{document}